\documentstyle[psfig,11pt,amsfonts,twoside]{article}
\newtheorem {theo} {\bf Theorem} [section]
\newtheorem {prop} [theo] {\bf Proposition}
\newtheorem {cory} [theo] {\bf Corollary}
\newtheorem {lem} [theo] {\bf Lemma}
\newtheorem {defn} [theo] {\bf Definition}

\newtheorem {exam} [theo] {\bf Example}
\newtheorem {rem} [theo] {\bf Remark}
\newcommand{\QED}{\hfill \CaixaPreta \vspace{6mm}}
\def\CaixaPreta{\vrule Depth0pt height6pt width6pt}

\newcommand{\qed}{\nopagebreak\hfill{\vrule width6pt height6pt depth0pt}}

\newcommand{\dpy}{\displaystyle}
\newcommand{\be}{\begin{eqnarray}}
\newcommand{\ee}{\end{eqnarray}}
\newcommand{\benn}{\begin{eqnarray*}}
\newcommand{\eenn}{\end{eqnarray*}}
\newcommand{\bse}{\begin{equation}}
\newcommand{\ese}{\end{equation}}
\newcommand{\bsenn}{\begin{displaymath}}
\newcommand{\esenn}{\end{displaymath}}

\newcommand{\logor}{\;\;{\rm or }\;\;}
\newcommand{\logif}{\;\;{\rm if }\;\;}
\newcommand{\logelse}{\;\;{\rm else }\;\;}
\newcommand{\then}{\;\;{\rm then }\;\;}
\newcommand{\where}{\;\;{\rm where }\;\;}
\newcommand{\with}{\;\;{\rm with }\;\;}
\newcommand{\such}{\;\;{\rm such\; that }\;\;}

\newcommand{\Hdim}{\;{\rm Hdim }\,}
\newcommand{\Hlim}{\;{\rm Hlim }\,}
\newcommand{\card}{\;{\rm card }\,}

\newcommand{\T}{\mbox{$\Bbb T$}} 	
\newcommand{\R}{\mbox{$\Bbb R$}} 	
\newcommand{\Z}{\mbox{$\Bbb Z$}} 	

\newcommand{\m}{mente }

\def\SBIMSMark#1#2#3{
 \font\SBF=cmss10 at 10 true pt
 \font\SBI=cmssi10 at 10 true pt
 \setbox0=\hbox{\SBF Stony Brook IMS Preprint \##1}
 \setbox2=\hbox to \wd0{\hfil \SBI #2}
 \setbox4=\hbox to \wd0{\hfil \SBI #3}
 \setbox6=\hbox to \wd0{\hss
             \vbox{\hsize=\wd0 \parskip=0pt \baselineskip=10 true pt
                   \copy0 \break%
                   \copy2 \break%
                   \copy4 \break}}
 \dimen0=\ht6   \advance\dimen0 by \vsize \advance\dimen0 by 8 true pt
                \advance\dimen0 by -\pagetotal
 \dimen2=\hsize \advance\dimen2 by .25 true in
  \openin2=publishd.tex
  \ifeof2\setbox0=\hbox to 0pt{}
  \else 
     \setbox0=\hbox to 3.1 true in{
                \vbox to \ht6{\hsize=3 true in \parskip=0pt  \noindent  
                {\it  Bol. Mex. Mat.~3}~{\bf 4} (1998) 1--24.
 
                \vfill}}
  \fi
  \closein2
  \ht0=0pt \dp0=0pt
 \ht6=0pt \dp6=0pt
 \setbox8=\vbox to \dimen0{\vfill \hbox to \dimen2{\copy0 \hss \copy6}}
 \ht8=0pt \dp8=0pt \wd8=0pt
 \copy8
 \message{*** Stony Brook IMS Preprint #1, #2 ***}
}

\newcommand{\multip}{\;\;{\rm multiplicity }\;\;}

\begin{document}
\title{HAUSDORFF DIMENSION OF BOUNDARIES OF SELF-AFFINE TILES IN $\R^N$}
\author{J. J. P. Veerman\\
Departamento de matem\'atica, UFPE, Recife, Brazil\thanks{e-mail:
veerman@dmat.ufpe.br}}
\date{\relax}
\maketitle
\SBIMSMark{1997/1}{January 1997}{}
\thispagestyle{empty}

\vskip .2in

\begin{abstract}
We present a new method to calculate the Hausdorff dimension of a certain class 
of fractals: boundaries of self-affine tiles. Among the interesting aspects are 
that even if the affine contraction underlying the iterated function system is 
not conjugated to a similarity we obtain an upper- and and lower-bound for its 
Hausdorff dimension. In fact, we obtain the exact value for the dimension if 
the moduli of the eigenvalues of the underlying affine contraction are all equal (this includes 
Jordan blocks). The tiles we discuss play an important role in the theory of wavelets. 
We calculate the dimension for a number of examples.
\end{abstract}

\vskip .3in

\section{Introduction}
\label{app=intro}
\setcounter{figure}{0}
\setcounter{equation}{0}

\vskip.2in
 
The object of this study is a class of self-affine (or self-similar) sets generated by pairs.
\begin{defn}
A pair $(M,R)$ is a linear isomorphism $M:\R^n \rightarrow \R^n$ with 
all eigenvalues outside the unit circle together with a finite subset $R$ of $\R^n$. 
\end{defn}
The space 
of closed and subsets of a fixed closed ball ${\cal B} \subset \R^n$ will be denoted by $H({\cal B})$. 
Endow this space with the usual Hausdorff distance between two compact sets 
(the infimum of $\epsilon$ such that an $\epsilon$-neighborhood of each one of the 
two sets contains the other). This distance induces a topology on $H({\cal B})$ 
with respect to which $H({\cal B})$ is a complete compact metric space. In $H({\cal B})$, 
we define 
\bsenn
\tau : H({\cal B}) \rightarrow H({\cal B})
\esenn
by
\bsenn
\tau(A) \buildrel\rm def \over \equiv \cup_{r\in R} M^{-1}(A+r) \quad .
\esenn
Such systems are affine examples of what are known as iterated function systems 
(see \cite{Ba}). It is easy to prove that $\tau$ is a contraction (see \cite{Hu}) 
and its unique fixed point is a compact set which we denote by $\Lambda(M,R)$. 
Hence we obtain that $\Lambda$ is `self-affine' (or `self-similar'):
\bse
\begin{array}{ll}
 &\Lambda = \cup_{r\in R} M^{-1}(\Lambda +r) \quad .\\[.4cm]
\logor &M\Lambda = \cup_{r\in R}(\Lambda +r) \quad . 
\end{array}
\label{eq=selfsimilar}
\ese

Equivalently (see \cite{HSV}), we may consider the set $\Lambda$ of expansions on the 
base $M$ using the set of digits $R$, or
\bsenn
\Lambda(M,R) = \{x \in \R^n | x = \sum_{i=1}^\infty M^{-i} r_i \with r_i \in R\} 
\quad .
\esenn 

Now we restrict our attention to a particularly interesting subclass of pairs, 
the class of standard pairs (in accordance with the 
nomenclature of \cite{LW1}). 
\begin{defn}
A standard pair $(M,R)$ is a pair such that the isomorphism $M$ 
preserves $\Z^n$ (that is: its matrix
has integer entries) and the set $R$ is contained in $\Z^n$ and is a 
complete set of coset representatives 
of $\Z^n/M\Z^n$ (R contains one representative
in $\Z^n$ of each the classes $\Z^n/M\Z^n$). 
\end{defn}
By performing a translation 
we may assume that $R$ contains the origin. Note that $R$ contains 
$m=|\det M|$ elements. This subclass of pairs is a boundary case in the 
following sense. In the above equation (\ref{eq=selfsimilar}), the measure 
on both sides can be positive, but only if the sets on the right hand side 
intersect in sets of (Lebesgue) measure zero.  

For later reference, we include the following definition. Denote by $\Z[M,R]$ the smallest $M$-invariant sublattice of $\Z^n$ that
contains all differences in $R$ (smallest in the sense that it contains no subset satisfying
the same requirements). 
\begin{defn}
A standard primitive pair is a standard pair with $\Z[M,R]=\Z^n$.
\end{defn}

\begin{defn}
Let $N:\R^n \rightarrow \R^n$ be a linear isomorphism preserving $\Z^n$ and 
$\pi_N : \R^n \rightarrow \R^n/N\Z^n$ the canonical projection. A compact 
set $A$ of positive Lebesgue measure is called a tile by $N\Z^n$ if $\pi_N : A \rightarrow \R^n/N\Z^n$ 
is a bijection for Lebesgue almost every point of $A$.
\end{defn}
When the matrix $N$ is not specified (as in most of this paper), we assume it 
to be the identity. In this case, we see that a tile is a compact set 
such that the union of its translates by $\Z^n$ covers $\R^n$, but any 
two translates by distinct elements of $\Z^n$ may intersect in sets of 
measure zero only.

The study of these tiles has important applications in various areas of mathematics. 
In fact, interest in these objects originally arose (see \cite{GM} and references 
therein), because they are intimately 
related to certain types of wavelets (Haar bases) (see \cite{D} for a discussion of 
wavelets). Connections with number theory (see \cite{LW4} and dynamical 
systems (see \cite{Ve2}) have been made by other authors. 

The main result here was proved by Lagarias and Wang. To state it, 
we need the following definition first.

\begin{defn}
A standard pair $(M,R)$ is called exceptional if there exists an integer matrix $P \in GL(n,\Z)$ such that
\begin{itemize}
\item $PMP^{-1}$ is a 'block-triangular' matrix $\left( \begin{array}{ll}
A & B \\ \emptyset & C \end{array} \right)$ and
\item $P(R)$ is of so-called quasi-product form (the definition is rather 
complicated, see \cite{LW5}).
\end{itemize}
A regular standard pair is one that is not exceptional.
\end{defn}

Here is the result of Lagarias and Wang \cite{LW5}.
\begin{theo}
\label{thm=lag-wang}
If $(M,R)$ is a regular standard pair then $\Lambda(M,R)$ is a tile by $\Z[M,R]$. 
If $(M,R)$ is an exceptional standard pair, then $\Lambda(M,R)$ is also a tile, 
though possibly by some other lattice. 
\end{theo} 

Lagarias and Wang based their proof upon an earlier work by Gr\"ochenig and Haas 
\cite{GH}, who proved a one dimensional version of this result. The proof by 
Gr\"ochenig and Haas was given a much more geometric flavor and was substantially 
simplified in \cite{Ve1}. 

The object of the current work is to study the Hausdorff dimension (for the 
definition, we refer to section 2 of this work) of the boundary 
$\delta \Lambda(M,R)$ of $\Lambda(M,R)$. In \cite{HSV}, it was proved that this boundary 
has Lebesgue measure zero. In the present 
work, we turn once again to the work by Gr\"ochenig and Haas. The calculation 
of dimension usually involves a counting argument. The counting procedure 
used in lemma \ref{thm=counting} is related to results of 
Gr\"ochenig and Haas (see \cite{GH}, lemma's 4.4 and 4.6). 

To state our main result formally, we need some notation. First, define the 
(direct) sum of two sets $A$ and $B$ as follows:
\bsenn
Z= A+B \buildrel\rm def \over \equiv \{z| z = a+b, a \in A, b \in B\} \quad. 
\esenn
Now define the difference set $D$:
\bsenn
D \buildrel\rm def \over \equiv R - R \buildrel\rm def \over \equiv 
\{d \in \Z^n|\exists r_1, r_2 \in R \such d = r_1 - r_2\} \quad ,
\esenn
and the multiplicity $\mu$:
\bsenn
\mu(d) = \card\{r_1,r_2 \in R| r_1 - r_2 = d\} \quad .
\esenn
(For instance, the element 0 occurs exactly $|\det M|$ times in $D$,
namely $0=r-r$ for all elements $r$ in $R$.) Denote 
\bse
S \buildrel\rm def \over \equiv \{\Lambda - \Lambda\} \cap \Z^n \quad ,
\label{eq=defS}
\ese
that is, the set of integer differences contained in $\Lambda$. Note that $S$ 
and $D$ are not equal! Let $\R^S$ be the space obtained by associating a fibre 
$\R$ to each point of $S$. Define a linear map $T : \R^S \rightarrow \R^S$, 
the transition operator, whose matrix elements are given by:
\bse
T_{ij} \buildrel\rm def \over \equiv \card\{(Mj+ D)\cap i)\} \quad ,
\label{eq=matrix}
\ese
where $i$ and $j$ are in $S$. Note that this matrix is the transpose of the one 
defined in \cite{GH} as contact-matrix and in \cite{Ve2} as transition operator.
The choice we make here is more natural in view of the calculations done in the 
next sections. 

Let us consider this transition operator in some more detail. First of all, 
$\frac{1}{m}T$ is a non-negative stochastic matrix:
\bse
\sum_{j\in S} T_{ij} = \sum_{j\in S} \card \{(Mj+D)\cap i\} = \card \{(M\Z^n+D)\cap i\} = m \quad .
\label{eq=leadingevalue}
\ese
The second equality follows from the definition of $S$: $i$ being in $S$ implies 
that if $j \not\in S$, then $\card \{(Mj+D)\cap i\} = 0$. The last equality 
is implied by the fact that $M\Z^n+D$ covers all elements of $\Z^n$ exactly $m$ 
times. 
Thus $T$ has leading eigenvalue $m$ with associated eigenspace $E_1 = \lambda(1,1,\cdots 1)$. 
(In fact, if $\Lambda$ is a tile, then all other eigenvalues have 
modulus smaller than $m$ see for example \cite{HSV}).

There is also a useful symmetry in the problem, namely, by its definition, $D$ is 
symmetric, that is $D = -D$. The same holds for $S$. It follows 
immediately from equation (\ref{eq=matrix}) that $T_{-i-j}=T_{ij}$. 
Thus $T$ preserves the subspace $E_+$ of symmetric 
vectors ($v_i=v_{-i}$) as well as the subspace of anti-symmetric vectors 
($v_i=-v_{-i}$) in $\R^S$. 
(Vectors in $\R^S$ are written as components $v_i$ where $i \in S$.) 

We decrease the dimensionality of the system by `quotienting out' this symmetry. 
Let $S^+ \subset S$ be such that of each pair $x\in S$ and $-x\in S$, 
precisely one is contained in $S^+$. Let $v^+$ denote the restriction 
of $v$ to $\R^{S^+}$. It is now easy to calculate the affine map 
$T^+ : \R^{S^+} \rightarrow \R^{S^+}$ induced by $T$ acting on $E_+$. Indeed, 
\benn
(T^+ v^+)_i &=& (\sum_{j\in S} T_{ij}v_j)^+ = \sum_{j\in S^+ -\{0\}}
(T_{ij}v_j + T_{i,-j}v_{-j}) + T_{i,\{0\}} v_{\{0\}} \\[.3cm]
&=& \sum_{j\in S^+ -\{0\}}(T_{ij} + T_{i,-j})v_j + T_{i,\{0\}} v_{\{0\}} 
\eenn
Notice that $T^+$ is a square matrix of dimension $\dpy\frac{|S|+1}{2}$ and one 
of its eigenvalues is $m$. The reader having trouble with these definitions can 
see them illustrated in an easy case (example \ref{example1}). 

The modulus of the eigenvalue of $M$ that is closest 
to the unit circle will denoted by $m_-$ (its reciprocal is the spectral radius 
of $M^{-1}$). Recall that $m = |\det M|$. Note that $1\leq m_-^n\leq m$.
\begin{defn}
\label{def=specialevalues} 
$\lambda$ is a special eigenvalue of $T$ if it is real, $\lambda$ is contained in 
$[m_-^{n-1},m)$, and $\lambda$ is an eigenvalue of $T^+$.
\end{defn}

Here is our main result.

\begin{theo}
\label{thm=main}
Let $(M,R)$ be a regular standard pair.\\
1) The Hausdorff dimension of $\delta \Lambda$ satisfies:
\bsenn
n + \frac{\ln \lambda - \ln m}{\ln m_-} \leq \Hdim(\delta \Lambda) 
\leq \frac{\ln \lambda}{ \ln m_-} \quad ,
\esenn
where $\lambda$ is the leading special eigenvalue of $T$ (which is the 
next-to-leading eigenvalue of $T$). \\
2) Let V be an open ball intersecting the boundary of $\Lambda$. 
The Hausdorff dimension of $\delta \Lambda\cap V$ is:  
\bsenn 
n + \frac{\ln \lambda_p - \ln m}{\ln m_-} \leq \Hdim(\delta \Lambda\cap V)  
\leq \frac{\ln \lambda_p}{ \ln m_-} \quad ,
\esenn 
where $\lambda_p$ is a special eigenvalue of $T$. 
\end{theo}

\noindent
{\bf Remark:}
In fact, the result also appears to hold for exceptional standard pairs, but 
we do not pursue this here.

Note that our result depends on the choice of $V$, that is: if $T$ has 
more than one positive real eigenvalue in the appropiate range, the dimension 
may not be constant. We will give an example of this in section \ref{examples} 
(example \ref{example8}). 

The most interesting case arises when the eigenvalues of $M$ are equal in 
modulus. We then have $m_-^{n} = m$, and the above inequalities become 
equalities.

\begin{cory}
Let $(M,R)$ be a regular standard pair and suppose that all eigenvalues of $M$ are 
equal in modulus.\\
1) The Hausdorff dimension of $\delta \Lambda$ satisfies:
\bsenn
\Hdim(\delta \Lambda) = \frac{\ln \lambda}{ \ln m_-} \quad ,
\esenn
where $\lambda$ is the leading special eigenvalue of $T$ (which is the 
next-to-leading eigenvalue of $T$).\\
2) For an open ball $V$ intersecting $\delta \Lambda$, we have 
\bsenn 
n-1 \leq \Hdim(\delta \Lambda\cap V) = \frac{\ln \lambda_p}{ \ln m_-} < n \quad ,
\esenn
where $\lambda_p$ is a special eigenvalue of $T$.
\end{cory}

In \cite{Ke}, Kenyon obtained an equality for the Hausdorff dimension similar 
to the one in the first part of the corollary. The difference is that in his case 
$\lambda$ is the principal eigenvalue of the transition operator of a 
Markov partition (in the usual sense) for $\delta \Lambda$. Our matrix $T$ is certainly not 
a transition matrix of a Markov partition, since it contains integers 
greater than one. Moreover, Kenyon's result does not include an algorithm to calculate  
the Markov transition operator from the initial data $(M,R)$. On the other hand, 
our transition operator is easily calculated (see section \ref{examples}, where 
we calculate the dimension of $\delta \Lambda$ in various cases). In addition, Kenyon's result is much more restricted than 
ours for various reasons. First of all, he assumes that $M$ is conformal, 
which in our corollary is not necessary (think of Jordan matrices). A more 
serious restriction is that he assumes that $\Lambda$ is homeomorphic 
to a ball in $\R^n$. This is generally not the case. A non-trivial set 
in one dimension, for example, cannot be connected. Thus his result 
excludes the one-dimensional tiles. In higher dimension, $\Lambda$ 
may have a complicated topology: in \cite{HSV} an example of a connected
set with infinitely many holes is given, or even one with 
infinitely many components, each of which has infinitely many holes
(see also example \ref{example8}). 

There is yet another approach to the calculation of the Hausdorff dimension for fractals 
generated by iterated function systems. This goes by means of a Mauldin-Williams 
graph (see \cite{MW}). This technique is also based on partitioning the boundary 
in a finite number of pieces and determining which ones are mapped where by the 
contractions of the iterated functions system. This sort of knowledge is not a priori present for the fractals under discussion here. In addition, the technique 
also requires the contractions to be conjugate to similarities.

Techniques by which one can calculate the Hausdorff dimension of a set that 
is invariant under a system of transformation that are not 
similarities are rare.  We know of only the socalled Sierpinski Carpets 
(see \cite{McM}, also explained in \cite{Fa2}). Here the transformations 
are {\it diagonal} matrices with integer entries. There is also an expression 
for the Hausdorff dimension of more general sets, sometimes called Falconer's Formula.
One can find this formula in \cite{Fa2}. 
This formula holds `almost always', but its proof does not indicate what 
the exceptional cases are (but see \cite{HL}).

Let $d_{\lambda_p}$ denote the size of the largest Jordan block associated with
$\lambda_p$. Denote the dimension of $\delta \Lambda \cap V$ by
\bse
\beta = \dpy\frac{\ln \lambda_p}{\ln m_-}
\label{eq=defbeta}
\ese
Denote the size of the largest Jordan block associated with $m_-$ by $d_M$.
Suppose $X$ is a set of Hausdorff dimension $\beta$. 
Denote by ${\cal H}^\beta (X)$ the Hausdorff outer measure of the set $X$ (for
the definition, see section 2). 

Our calculations needed for the above result also give the following result.

\begin{theo}
\label{thm=ssets}
Let $(M,R)$ be a regular standard primitive pair and suppose all
eigenvalues of $M$ have equal modulus. \\
1) If $d_M= d_{\lambda_p}=1$, then:
\bsenn
{\cal H}^\beta (\delta\Lambda \cap V) < \infty \quad .
\esenn
2) If $d_{\lambda_p} -1 \geq (n-\beta)(d_M -1)$, then
\bsenn
{\cal H}^\beta(\delta\Lambda \cap V) > 0 \quad .
\esenn
3) If $d_{\lambda_p} -1 > (n-\beta)(d_M -1)$, then
\bsenn
{\cal H}^\beta(\delta\Lambda \cap V) = \infty \quad .
\esenn
\end{theo}

The set-up of this article is as follows. In section 2 we discuss elementary 
notions concerning the set $\delta \Lambda$. Essentially, we decribe how to 
construct it. Then in the next section, we describe the operator $T$ and 
its properties. This serves to facilitate the counting argument already mentioned. 
These counting arguments will then be spelled out in the next two sections. The
first of these, section 4, gives the upper bound for our dimension estimate, and 
in section 5 we derive the lower bound. In section 6, we prove the result concerning 
the regularity of the boundary (if $M$ has eigenvalues of equal modulus). Finally, 
in section 7, we calculate the dimension of the boundary of several tiles.

In all subsequent sections, we will assume, without loss of generality, that $(M,R)$ is a standard 
 primitive pair. We give the reduction to that case here.

\begin{lem}
Let $(M,R)$ be a standard pair. Then $\Lambda(M,R)$ is conjugate to a tile 
$\Lambda'=\Lambda(M',R')$, where $(M',R')$ is a standard primitive pair.
\end{lem}

\noindent
{\bf Proof:} Let $B$ be a matrix such that $B\Z^n=\Z[M,R]$. Clearly, $M$ 
preserves $B\Z^n$. Thus $M'=B^{-1}MB$ has integer entries. Further, $R'=B^{-1}R 
\subset \Z^n$ is a complete set of coset representatives \cite{Ve1}. It is easy 
to see that 
$\Lambda'=B^{-1}\Lambda$. 
\QED

The lattice $B\Z^n$ is a generalization of the notion of greatest common divisor 
(see also \cite{Ve1}).

\vskip .2in
\noindent
{\bf Acknowledgements:} I am grateful to Leo Jonker for reading preliminary 
versions of this man\-u\-script and suggesting many improvements. I also thank 
Jos\'e Ruidival dos Santos for bringing 
proposition \ref{thm=alaoglu} to my attention and to Geraldo Oliveira and Andr\'e Banks Rocha 
for useful conversations. 

\pagebreak 

\section{Construction of the Boundary}
\label{construction}
\setcounter{figure}{0}
\setcounter{equation}{0}
 
\vskip.2in

In this section, we describe the construction of the boundary. It also 
serves to collect most of the definitions and notation (in so far 
not already discussed in the introduction) in addition 
to some elementary results. 

Throughout this work, we will denote Lebesgue measure by $\mu$ and 
the $\epsilon$-neighborhood of a set $A$ by $N_\epsilon (A)$.
The diameter of a set $U$ is denoted by $|U|$.

Of the many definitions of dimension \cite{Ma}, Hausdorff dimension has been 
one of the mathematically most fruitful ones (see \cite{Fa2}). We give the 
definition and refer to \cite{Fa1} and \cite{Fa2} and references therein for 
further reading. 

For a set $F\subset \R^n$, define
\bsenn
{\cal H}^s_\delta (F) = \inf \{\sum |U_i|^s \} \quad ,
\esenn
where the infimum is over all countable covers of $F$ whose individual sets have 
diameters less than $\delta$. The s-dimensional Hausdorff outer-measure is given by:
\bsenn
{\cal H}^s (F) = \lim _{\delta \rightarrow 0} {\cal H}^s_\delta (F) \quad .
\esenn
There is a unique number $\beta \geq 0$ such that
\bsenn
\begin{array}{rcr}
\logif s < \beta & \then & {\cal H}^s (F) = \infty \\[.3cm]
\logif s > \beta & \then & {\cal H}^s (F) = 0
\end{array} \quad .
\esenn
This number is the Hausdorff dimension of $F$. Part of the importance derives 
from the fact that it is an invariant under bi-Lipschitz homeomorphisms. A set is called 
an $s$-set if its $s$-dimensional Hausdorff outer-measure is positive and finite.
This property is also invariant under bi-Lipschitz homeomorphisms.

\begin{lem}
\label{thm=hasmeasure}
For all $x \in \Lambda$, we have that for all $\epsilon > 0$
\bsenn
\mu(N_\epsilon(x)\cap \Lambda) > 0 \quad .
\esenn
\end{lem}

\noindent
{\bf Proof:} By definition of $\Lambda$, for any given $\epsilon$, 
there is a $k$ such that $x_k+M^{-k}\Lambda
\subset N_\epsilon(x)$ for some $x_k \in \sum_{i=1}^k M^{-i} R$.
\QED
 
\begin{cory}
$\Lambda$ is the closure of its interior.
\end{cory}

\noindent
{\bf Proof:} Recall that $\delta \Lambda$ has measure zero and use the previous
lemma.
\QED

Now we define 
\bse
\Lambda^{(2)} = \{ x \in \Lambda | x + \{\Z^n - \{0\}\} \cap \Lambda \neq \emptyset
\} \quad .
\label{eq=lambda2}
\ese
Recall that we assume that $(M,R)$ is a regular standard primitive pair. Thus 
by theorem \ref{thm=lag-wang}, $\Lambda(M,R)$ is a tile by $\Z^n$.

\begin{lem} 
\label{thm=lam2=deltalam}
We have
\bsenn
\Lambda^{(2)} = \delta \Lambda \quad .
\esenn
\end{lem}

\noindent
{\bf Proof:} First, let $x \in \delta \Lambda$. Then, for all $\epsilon >0$
\benn
\pi[\{N_\epsilon(x)+\Z^n\}\cap \Lambda] = \pi N_\epsilon(x)\\[.4cm]
\pi[N_\epsilon(x)\cap \Lambda] \neq \pi N_\epsilon(x) \quad .
\eenn
Thus for all $\epsilon > 0$, $[N_\epsilon(x)+\{\Z^n-\{0\}\} \cap \Lambda$ is non-empty.

Now let $x\in \Lambda^{(2)}$. There is $y\in x + \{\Z^n-\{0\}\}$ in $\Lambda$.
Pick some small $\epsilon_0$ such that
\bsenn
N_{\epsilon_0}(x)\cap N_{\epsilon_0}(y) = \emptyset
\esenn
For all $\epsilon<\epsilon_0$ we then have
\bsenn
\pi[(N_\epsilon(x)\cap \Lambda) \cup (N_\epsilon(y)\cap \Lambda)]
\esenn
covers $\pi[N_\epsilon(x)]$ at most once (modulo sets of measure zero). 
By the previous lemma,
$N_\epsilon(x)\cap \Lambda$ and $N_\epsilon(y)\cap \Lambda$
have positive measure. Thus, for all $\epsilon > 0$, neither has full
measure.
\QED
 
We define: 
\bse
\label{eq=defgamma}
\Gamma_k \buildrel\rm def \over \equiv \sum_{i=1}^{k}M^{-i}R = \tau^k(\{0\}) \quad .
\ese
$\Gamma_k$ is the $k$-th generation approximation to $\Lambda$ and as an element 
of $H({\cal B})$ it is the $k$-th iterate of $\{0\}$. 
Notice that since $0\in R$, $\Gamma_k \subseteq \Lambda$. 
Similarly, recall the definition 
of the set $S$ (equation (\ref{eq=defS})) and define 
\bse
\label{eq=defdelta}
\Delta_k \buildrel\rm def \over \equiv \{x\in \Gamma_k|\exists i\in S-\{0\}, \exists j\in S-\{0\} \such x+i+M^{-k}j\in \Gamma_k\} \quad ,
\ese
the $k$-th generation approximation to the boundary of $\Lambda$. $E_k$ is the 
$k$-th generation `filled-in' approximation to the boundary of $\Lambda$: 
\bse
\label{eq=defE}
E_k \buildrel\rm def \over \equiv \Delta_k + M^{-k}\Lambda \quad .
\ese

\begin{lem}
\label{thm=inclusion1}
We have 
\bsenn
\Delta_k \subseteq N_{3\epsilon_k}(\delta\Lambda)
\esenn
where $\epsilon_k=|M^{-k}\Lambda|$
\end{lem}

\noindent
{\bf Proof:} If $x\in \Delta_k$, then $y=x+i+M^{-k}j\in \Delta_k$. So, both 
$x$ and $y$ are in $\Lambda$. Then $N_{\epsilon_k}(y)$ contains an open set 
$A \subset \Lambda$ lemma \ref{thm=hasmeasure}. The set $A-i\subset N_{3\epsilon_k}(x)$ is not contained 
in $\Lambda$. Since $x\in \Lambda$, $N_{3\epsilon_k}(x)$ must also contain points 
of the boundary. 
\QED

We now construct a sequence of maps $\tau_k : H({\cal B})\rightarrow H({\cal B})$ whose 
fixed points $F_k$ will also approach the boundary of $\Lambda$:
\bse
\label{eq=deftauk}
\tau_k(X) \buildrel\rm def \over \equiv M^{-k}X + \Delta_k \quad ,
\ese
and 
\bse
\label{eq=defF}
F_k \buildrel\rm def \over \equiv \Lambda(M^k, M^k \Delta_k) \quad .
\ese
Note that $\tau_k F_k = F_k$.

\begin{lem} We have
\label{thm=inclusion2}
\bsenn
\delta \Lambda \subseteq F_k \subseteq E_k \quad .
\esenn
\end{lem}
 
\noindent 
{\bf Proof:} The second inclusion follows immediately from the definitions
of $E_k$ and $F_k$. 
 
Suppose $x_1\in \delta\Lambda$, then by lemma \ref{thm=lam2=deltalam}, 
we may suppose that $x_1 \in \Lambda^{(2)}$. We will
show that $x_1$ can be written as
\bsenn
x_1 = \sum M^{-ki}v_i , \where v_i \in M^k \Delta_k \quad .
\esenn
By iterating equation (\ref{eq=selfsimilar}), we see that there is a $v_1 \in M^k \Gamma_k$ 
with
\bsenn
x_1 \in M^{-k}(\Lambda + v_1) \quad .  
\esenn
Thus there is an $x_2 \in \Lambda$ such that 
\bsenn
x_1 = M^{-k}(x_2 + v_1) \quad .
\esenn
Since $x_1 \in \Lambda^{(2)}$, there are ${\overline x}_1 \in \Lambda$
and ${\overline v}_1 \in M^k \Gamma_k$ such that for some 
${\overline x}_2 \in \Lambda$
\bsenn
{\overline x}_1 = M^{-k}({\overline x}_2 + {\overline v}_1) \quad ,\quad 
x_1 - {\overline x}_1 = j_1 \in S -\{0\} \quad .
\esenn
Now observe that 
\bsenn
v_1 - {\overline v}_1 = M^k(x_1 - {\overline x}_1) + {\overline x}_2 - x_2 \quad .
\esenn
Since $v_1 - {\overline v}_1 \in \Z^n$, we also have
\bsenn
v_1 - {\overline v}_1 = M^k j_1 + j_2 \quad ,
\esenn
Note that $v_1 - {\overline v}_1 \in M^k(\Gamma_k - \Gamma_k)$. This set contains
no elements of $M^k(\Z^n - \{0\})$ (see \cite{HSV}). Recall that $j_1 \in S-\{0\}$. 
Thus $j_2 \in S-\{0\}$. 
Thus $v_1 \in M^k\Delta_k$. Moreover, we see that
\bsenn
{\overline x}_2 - x_2 = j_2 \quad .
\esenn
Thus $x_2$ belongs to $\delta \Lambda = \Lambda^{(2)}$. 
Continue by induction.
\QED
 
The following proposition establishes that the sets we have defined so far
converge to $\delta \Lambda$ in the Hausdorff topology. 

\begin{prop}
\label{thm=boundary}
One has 
\bsenn
\Hlim_{k\rightarrow \infty} F_k = \Hlim_{k\rightarrow \infty} \Delta_k = \Hlim_{k\rightarrow \infty} E_k =
\delta \Lambda \quad .
\esenn
\end{prop}
 
\noindent 
{\bf Proof:} Combining lemmas \ref{thm=inclusion1} and \ref{thm=inclusion2}, we have 
\bsenn
\Delta_k \subseteq N_{3\epsilon_k}(\delta \Lambda) \subseteq N_{3\epsilon_k}(F_k) 
\subseteq N_{3\epsilon_k}(E_k) \quad .
\esenn
By definition of the Hausdorff distance, the Hausdorff distance between 
the two sets $\Delta_k$ and $E_k$ is exactly
$|M^{-k}\Lambda|$, which implies the proposition.
\QED

\pagebreak

\section{The Transition Matrix}
\label{transition}
\setcounter{figure}{0}
\setcounter{equation}{0}

\vskip .2in

Here we derive how the transition operator counts the number of points 
in the $k$-th level approximation of the boundary. Again, we assume that 
$\Lambda$ is a tile (for example, when $(M,R)$ is a regular standard primitive 
pair, according to theorem \ref{thm=lag-wang}). 
This implies that the leading eigenvalue of $T$, discussed
in equation (\ref{eq=leadingevalue}) is simple (see \cite{GH} or \cite{Ve1}). 
 
Consider the transition matrix $T : \R^S \rightarrow \R^S$ as defined in the 
introduction. For an open set $V$ with diameter smaller than one and a non-negative integer $k$ define the 
contact-matrix (we borrowed the name from \cite{GH}) as follows.
\bse
T(k,V)_{ij} \buildrel \rm def \over \equiv \left\{
\begin{array}{ll}
\card\{x\in \Delta_k | x \in V \logor x+i+M^{-k}j\in V\} & \logif i,j \in S-\{0\}\\
0 & \logelse
\end{array} \right.
\label{eq=defcontact} 
\ese
In this definition, we say that $x$ is a basepoint for the difference 
$i+M^{-k}j$. 

Note that we have approximately
\bsenn
T(k+\ell,V) \approx T^\ell T(k,V) \quad .
\esenn
The fact that this is not exact, is due to boundary effects whose relative error 
decreases exponentially (this is worked out precisely in proposition \ref{thm=boundaryeffect}). 
By applying $T$ to $T(k,V)$, one sees that 
\bsenn
(T^\ell T(k,V))_{\{0\},j}= (T^\ell T(k,V))_{i,\{0\}}= 0 \quad .
\esenn
Thus for all $\ell$, the angle between the span of the columns of $T^\ell T(k,V)$ and  
$E_1= \lambda (1,1\cdots ,1)$, the eigenspace associated with the leading eigenvalue $m$, is bounded from 
below. Since, as remarked just after equation (\ref{eq=leadingevalue}), the leading eigenvalue 
is simple, the growth-rate of $T^\ell T(k,V)$ is less than $m$. 

Since we will be doing a lot of counting, define the following counter for a 
non-negative integer matrix C:  
\bsenn
\|C\| = \sum_{ij} c_{ij} \quad .
\esenn

For a given ball $B$, we can now express the growth-rate of 
$\card(\Delta_{k+\ell}\cap B)$ in terms of the growth-rate associated with 
the matrix $T$. This is done in the following two results.

\begin{lem} We have
\label{thm=counting}
\bsenn
\frac{\|T(k,B)\|}{2 \card S} \leq \card(\Delta_k\cap B) \leq \frac{\|T(k,B)\|}{2}
\quad .
\esenn
\end{lem}
 
\noindent
{\bf Proof:} The first inequality follows from the fact $x \in \Delta_k\cap B$ can 
be the basepoint of at most $2 \card S$ differences.

The second inequality follows from the definition of $\Delta_k$:
each $x\in \Delta_k\cap B$ is the basepoint of some difference and its negative (by the definition 
of $\Delta_k$). 
\QED
 
\begin{prop}
\label{thm=boundaryeffect} 
Let $B_r$ be a ball of radius $r$ intersecting $\delta \Lambda$. Fix a 
constant $\epsilon > 0$ and choose $k$ such that $|M^{-k}\Lambda| < \epsilon r$.
Then
\bsenn
(T^{\ell}T(k,B_{r(1-\epsilon)}))_{i,j} \leq (T(k+\ell,B_r))_{i,j} \leq (T^\ell T(k,B_{r(1+\epsilon)})_{i,j} \quad .
\esenn
\end{prop}
 
\noindent
{\bf Proof:} Suppose $\Delta_k\cap B$ contains $N = T(k,B)_{ba}$ basepoints of the 
difference $a+M^{-k}b$. Then $\Delta_k\cap B + \sum_{i=k+1}^{k+\ell} M^{-i}R$
contains
\benn
\sum_{b\in S} \card \left( \{a+M^{-k}b + \sum_{i=k+1}^{k+\ell}M^{-i}R\} \cap 
\{a+M^{-k-\ell}c +\sum_{i=k+1}^{k+\ell}M^{-i}R\}\right) \cdot T(k,B)_{ba} = \\[.4cm]
=\sum_{b\in S} \card(\{b + \sum_{i=1}^{\ell} M^{-i}D\} \cap \{M^{-\ell}c\}) \cdot  
T(k,B)_{ba} =\sum_{b\in S} (T^\ell)_{cb} T(k,B_r)_{ba}
\eenn
basepoints of the  difference $a+M^{-k-\ell}c$. There is a discrepancy due to the 
fact that $\Delta_k\cap B + \sum_{i=k+1}^{k+\ell} M^{-i}R \neq \Delta_{k+\ell}\cap B$.
Points may `seep' across the boundary of $B$. By assumption, points that do so,
lie within a distance $\epsilon r$ of the boundary.
\QED
 
\begin{lem}
\label{thm=specialevalues}
Let $V$ be a ball.\\
1) the growth-rate (in $\ell$) of $\|T^\ell T(k,V)\|$ is 
determined by a special eigenvalue. \\
2) If $V$ is sufficiently large, then this eigenvalue is the next-to-leading 
eigenvalue of $T$.
\end{lem}

\noindent
{\bf Proof:} To prove the first statement, recall the definition of the 
$T$-invariant (symmetric and anti-symmetric) splitting 
of $\R^n$ in the introduction: $\R^n = E_+ \oplus E_-$, and the 
operators $T_+$ and $T_-$ which are just the linear map $T$ restricted
to these respective spaces. Let $v=v^+ + v^-$ be the $j$-th column of 
$T(k,V)$. The $j$-th column of $T^\ell T(k,V)$ is
\bsenn
(T^+)^\ell v^+ + (T^-)^\ell v^- \quad .
\esenn
By the previous proposition, one sees that the components of this
vector are non-negative for all $\ell$. Thus the growth-rate must be 
determined by an eigenvalue of $T^+$.

A similar argument shows that this eigenvalue is real positive. Let $E_i$ denote
the eigenspaces of $T^+$, ordered in such a way that the associated
eigenvalues $\lambda_i$ satisfy: $|\lambda_i| \geq |\lambda_{i+1}|$. 
Denote $T|_{E_i}$ by $T_i$. Now, consider the smallest integer $j$ 
for which the span of the columns of $T^\ell T(k,V)$ intersects $E_j$
in a linear subspace of positive dimension. The growth of some 
column $v = v_j + \sum_{\alpha > j} v_\alpha$ where $v_\alpha \in E_\alpha$,
is dominated by
$(T_j)^\ell v_j$. Suppose that the eigenvalues $\lambda_{\alpha>j}$ 
are less in modulus than $\lambda_j$. Then for $\ell$ big enough the
entries of $(T_j)^\ell v_j$ have to be positive, since its contribution 
dominates the count of the number of differences.

Finally, the bounds on the eigenvalue follow from theorem \ref{thm=main} 
(and will not be used in the proof of that theorem), together with 
the observation that the boundary of an $n$-dimensional volume 
has dimension at least $n-1$. 

To prove the second statement, consider the following splitting:
\bsenn
\R^n = E_1 \oplus E_\perp \quad .
\esenn
If $V$ sufficiently big, it will contain $\Lambda$. In this case, 
$T(k,V)$ counts all differences in $\Gamma_k$ except the ones 
equal to zero. Thus 
for any vector $v$, we have that
\bsenn
T(k,V)v = T(k,V)(v_1\oplus v_\perp) = T^k v_\perp \quad .
\esenn
\QED
 
Finally, we state a lemma that we will often use

\begin{lem}
\label{thm=growthrate}
Let $A:\R^n \rightarrow \R^n$ a linear map and $E_{\lambda}$ the invariant space 
associated with an eigenvalue of modulus $\lambda$. \\
1) There is a polynomial $p$
such that for all $x\in E_{\lambda}$ 
\bsenn
C_1 \lambda^k |x| \leq |A^kx| \leq C_2 \lambda^k p(k) |x| \quad .
\esenn
The degree of $p$ is one less than the size of the biggest Jordan block 
associated with $A$.\\
2) For a given $x \in E_{\lambda}$, we have that there is a polynomial $p$ of degree 
less than the size of the biggest Jordan block associated with $A$ such that
\bsenn
C_1 \lambda^k p(k) \leq |A^kx| \leq C_2 \lambda^k p(k) \quad .
\esenn
\end{lem}

\noindent
{\bf Proof:} Bring $A$ into Jordan form.
\QED

\pagebreak

\section{The Upper Bound for the Dimension}
\label{upper}
\setcounter{figure}{0}
\setcounter{equation}{0}
 
\vskip.2in

In this section, we calculate the upper bound of the dimension of $\delta \Lambda 
\cap V$ where $V$ is an arbitrary open ball intersecting the boundary of a tile 
$\Lambda$. We do this by showing there is a sequence of sets $E_k$ with 
$\delta \Lambda\cap V \subseteq E_k$. The upper bound we calculate equals
$\lim_{k\rightarrow \infty} \Hdim E_k$. 

The main tool we use to give an upper estimate for the dimension of a set 
is one that follows almost directly from proposition 9.6 in \cite{Fa2}. For 
completeness we include the proof. 

Notice that by lemma \ref{thm=growthrate}, there is a polynomial $p$ of degree 
$d_{m_-}$ such that $|M^{-k}\Lambda| \leq C p(k) m_-^{-k}$. 

\begin{prop}
\label{thm=upper}
Let $(A,Q)$ be a pair as described in the introduction (not necessarily 
standard) and denote the spectral radius of $A^{-1}$ by $a_-^{-1}$. Then
\bsenn 
\Hdim \Lambda(A,Q)  \leq \frac{\ln \card (Q)}{\ln a_-} \quad .
\esenn
\end{prop}

\noindent
{\bf Proof:} Cover $\Lambda(A,Q)$ by hypercubes whose sides have length 
$\epsilon_k = |A^{-k}\Lambda| \leq C p(k) a_-^{-k}$.  
To do so, we need at most $(\card Q)^k$ hypercubes. Let $\beta = \dpy\frac{\ln \card (Q)}{\ln a_-}$ 
and $d$ any positive number. Then the $\beta+d$ dimensional Hausdorff outer measure 
of $\Lambda$ satisfies:
\benn
{\cal H}_{\epsilon_k}^{\beta+d}(\Lambda) &\leq& C p(k)^{\beta+d} a_-^{-k(\beta+d)} (\card Q)^k \\[.4cm]
& = & C p(k)^{\beta+d} a_-^{-kd} \quad ,
\eenn
which tends to zero as $k$ tends to infinity (and $\epsilon_k$ to zero). 
Thus for all positive $d$, 
the $\beta+d$-dimensional Hausdorff outer measure of $\Lambda$ is zero. The 
proposition follows immediately from the definition of the Hausdorff dimension. 
\QED

\begin{cory}
\label{thm=upper2}
Let $(M,R)$ be a regular standard primitive pair and 
suppose that $Q$ is a subset of $\cup_{i=0}^{k-1} M^i R$. Then 
\bsenn
\Hdim \Lambda(M^k,Q) \leq \frac{\ln \card (Q)}{k \ln m_-} \quad .
\esenn 
\end{cory}

\noindent
{\bf Remark:} 
This estimate becomes exact if $M$ is conformal and 
the pair $(M,R)$ satisfies 
the open set condition (see \cite{Fa1}). As remarked in the introduction, 
$\Lambda(M,R)$ is a tile whose 
boundary has measure zero. Thus its interior has positive measure. The open 
set condition now holds for the pair $(M^k,Q)$ with the interior of $\Lambda$ 
as the open set.

Notice that by lemma \ref{thm=growthrate} there is a polynomial $q$ of degree 
$d_{\lambda_p}-1$ such that 
\bsenn
\|T^\ell T(k,V)\| \leq C q(\ell) \lambda_p^\ell \quad .
\esenn
 
\begin{lem}
\label{thm=countdelta}
There is a polynomial $q$ such that 
\bsenn
C_1 q(k) \lambda_p^k < \card(\Delta_k\cap V) \leq C_2 q(k) \lambda_p^k \quad ,
\esenn
where $\lambda_p$ is a special eigenvalue (next-to-leading if $V$ big enough).
\end{lem}

\noindent
{\bf Proof:} Let $V_+$ denote the $\epsilon r$-neighborhood of $V$ and $V_-$ the ball
whose $\epsilon r$-neighborhood is $V$. Then by results
\ref{thm=counting}, \ref{thm=specialevalues} and \ref{thm=boundaryeffect}
\bsenn
\|T^\ell\cdot T(k,V_-)\| \leq \card(\Delta_{k+\ell}\cap V) \leq \|T^\ell\cdot T(k,V_+)\| \quad .
\esenn
The result follows from applying lemmas \ref{thm=specialevalues} and \ref{thm=growthrate} to this formula. 
\QED

\begin{theo} 
$\Hdim(\delta \Lambda \cap V) \leq \dpy\frac{\ln \lambda_p}{\ln m_-} \quad .$
\end{theo}

\noindent 
{\bf Proof:} We have by lemma \ref{thm=inclusion2}
\bsenn
\delta \Lambda \cap V \subseteq F_k\cap V \quad .
\esenn
Thus,
\bsenn
\Hdim(\delta \Lambda\cap V) \leq \liminf_{k\rightarrow \infty} \Hdim (F_k\cap V) \quad .
\esenn
Now apply corollary \ref{thm=upper2} and lemma \ref{thm=countdelta} to the definition 
of $F_k$ (equation (\ref{eq=defF})) to see that 
\bsenn
\liminf_{k\rightarrow \infty} \Hdim F_k \leq \liminf_{k\rightarrow \infty} \dpy\frac{\card(\Delta_k\cap V)}{k\ln m_-} = 
\dpy\frac{\ln \lambda_p}{\ln m_-} \quad .
\esenn
\QED

\pagebreak

\section{The Lower Bound for the Dimension}
\label{lower}
\setcounter{figure}{0}
\setcounter{equation}{0}
 
\vskip.2in

In this section, we calculate the lower bound for the Hausdorff dimension of 
$\delta \Lambda\cap V$, where $V$ is an open ball intersecting the boundary. 

The technique we use consists putting a probability measure $\nu$ on 
$\delta \Lambda\cap V$ where $V$ is an open ball. If for $r$ small enough
the measure contained in a ball of radius $r$ does not exceed $Cr^s$, then 
$s$ is a lower bound for the Hausdorff dimension \cite{Fa2}. The following result 
is a minor extension of this.

\begin{prop}
\label{thm=falconer}
Let $\nu$ be a probability measure on a set $X$. If for all $d>0$ there 
exists a $C_d >0$ such that for all $x\in X$ 
\bsenn
\lim_{r\rightarrow 0} \dpy\frac{\nu(B_r(x))}{r^{\beta-d}} \leq C_d \quad .
\esenn
Then
\bsenn
{\cal H}^{\beta-d}(X) \geq \dpy\frac{\nu(X)}{C_d} \quad .
\esenn
And thus $\Hdim (X) \geq \beta$.
\end{prop}

\noindent
{\bf Proof:} Let $\{U_i\}$ be a cover of $X$ and suppose that $x_i\in U_i\cap X$. 
Then
\bsenn
\nu(U_i) \leq \nu(B_{|U_i|}(x_i)) < (C_d+\epsilon) |U_i|^{\beta-d} \quad ,
\esenn
provided $|U_i| < \delta$ and $\delta$ small enough. By summing, one obtains:
\bsenn
\sum_{i} |U_i|^{\beta-d} > \dpy\frac{\nu(X)}{C_d + \epsilon} \quad .
\esenn
The definition of Hausdorff dimension implies that for all $d>0$
\bsenn
\Hdim(X) \geq \beta - d \quad .
\esenn
\QED

We now define a probability measure $\nu_k$ on $\Delta_k\cap V$. For any open
set $U \subset V$:
\bse
\nu_k(U) = \frac{\card (\Delta_k\cap U)}{\card (\Delta_k\cap V)} \quad .
\label{eq=measure}
\ese
Presumably, the measures $\nu_k$ converge exponentially fast to a limiting 
measure $\nu$. However, the local structure of $\Delta_k$ is difficult to control. 
Instead, we simply use the Banach-Alaoglu theorem (see \cite{Fr}). 

\begin{prop}
\label{thm=alaoglu}
There is a subsequence $\{\nu_{k_i}\}$ that converges to a probability measure
$\nu$ on $\delta \Lambda\cap V$.
\end{prop}

\noindent
{\bf Proof:} Let $X= \cup_k (\Delta_k \cap V)$. Then $X$ is compact and $\nu_k(X)=1$. 
By the Banach-Alaoglu
theorem, there must be a subsequence of the $\nu_k$ converging to a measure $\nu$. By using 
proposition \ref{thm=boundary}, we see that $\nu$ has support in $\delta \Lambda\cap V$.
\QED

We will now use proposition \ref{thm=boundaryeffect} to devise a method to count the growth rate with $k$ of the number 
of elements of $\Delta_k$ contained in a ball $B$. This will enable us to 
calculate estimates for the measures $\nu_k$ and thus to determine the $\nu$-measure 
contained in a ball of radius $r$. Without loss of generality, we may also take $V$ to be a ball. 

\begin{prop}
Let $\nu$ be the probability measure just constructed. Fix a small constant
$\epsilon>0$. Then there is a $K>0$ such that if $k$ satisfies
\bsenn
|M^{-k}\Lambda| < \epsilon r \quad ,
\esenn
then for a ball of radius $r$ intersecting $\delta \Lambda \cap V$ 
\bsenn
\nu(B_r) < K r^n \frac{m^k}{\lambda_p^k} \frac{1}{q(k)} \quad ,
\esenn
where $\lambda_p$ is a special eigenvalue (next-to-leading if $V$ big enough). 
\end{prop}

\noindent
{\bf Proof:} Let $V_-$ be the ball whose $\epsilon r$-neighborhood equals $V$, 
denote the $\epsilon r$-neighborhood of $B_r(x)$ by $B_+$, and assume that 
$B_+ \subset V_-$. Combining  the hypotheses 
and the results \ref{thm=counting}, \ref{thm=specialevalues}, and \ref{thm=boundaryeffect}, we see
\bsenn
\card(\Delta_{k+\ell} \cap V) \geq \dpy\frac{\|T^\ell \cdot T(k,V_-)\|}{\card S} \quad .
\esenn
Using the opposite inequalities, we arrive at 
\bsenn
\card(\Delta_{k+\ell} \cap B_r) \leq \|T^\ell \cdot T(k,B_+)\| \quad .
\esenn
Now, since $(\Delta_{k+\ell}\cap B_+)\subset (\Delta_{k+\ell}\cap V_-)$, the 
growth rate (in $\ell$) of 
$T^\ell \cdot T(k,B_+)$ is dominated by the growth rate of 
$T^\ell \cdot T(k,V_-)$. 
Thus, recalling the definition of the measure $\nu$ in proposition \ref{thm=alaoglu}, 
and lemma \ref{thm=counting}:
\bsenn
\nu(B_r) \leq \lim_{i} \frac{\card(\Delta_{k+\ell_i} \cap B_r)}{\card(\Delta_{k+\ell_i} \cap V)} 
\leq \card S\frac{\|T(k,B_+)\|}{\|T(k,V_-)\|} 
\esenn
By lemma \ref{thm=countdelta}, we know that 
\bsenn
\|T(k,V_-)\| \geq C \lambda_p^k q(k) \quad .
\esenn
Furthermore,
\bsenn
\|T(k,B_+)\| \leq \dpy\frac{{\rm vol}(B_+)}{{\rm vol}(M^{-k}\Lambda)}  = \frac{C r^n}{m^{-k}} \quad .
\esenn
Putting the estimates together yields the result.
\QED

To simplify notation, put
\bse
\beta = n + \dpy\frac{\ln \lambda_p - \ln m}{\ln m_-} \quad .
\label{eq=beta}
\ese

\begin{theo}
\label{thm=lowerbound}
$\Hdim (\delta \Lambda \cap V) \geq \beta $.
\end{theo}

\noindent
{\bf Proof:} In accordance with the previous proposition, we can choose $k$ such that 
\bse
|M^{-k}\Lambda| < \epsilon r \leq |M^{-(k-1)}\Lambda| \quad .
\label{eq=choicek}
\ese
By bringing $M^{-1}$ in Jordan normal form applying lemma \ref{thm=growthrate}, 
we conclude that there are a constant $C$ and a polynomial $p$ of degree $d_M - 1$ such that 
$C_1 p(k) m_-^{-k} < |M^{-(k-1)}\Lambda| \leq C_2 p(k) m_-^{-k}$. By using this and equation \ref{eq=choicek}, 
we calculate the dependency of $k$ on $r$: 
\bse
C_1 \epsilon^{-1} p(k) m_-^{-k} < r \leq C_2 \epsilon^{-1} p(k) m_-^{-k} \quad .
\label{eq=findr}
\ese
By proposition \ref{thm=falconer}, we will be done if for all positive $d$ 
small enough
\bsenn
\lim_{r\rightarrow 0} \dpy\frac{\nu(B_r(x))}{r^{\beta-d}} \leq C_d \quad .
\esenn
We calculate, using the previous proposition and equation (\ref{eq=beta}):
\benn
\dpy\frac{\nu(B_r(x))}{r^{\beta-d}} < Kr^n m^k \lambda_p^{-k} {\frac{1}{q(k)}} r^{-n} r^{-\frac{\ln \lambda_p}{\ln m_-}} r^{\frac{\ln m}{\ln m_-}} r^d = \\[.3cm]
K m^k \lambda_p^{-k} m^{\frac{\ln r}{\ln m_-}} \lambda_p^{\frac{-\ln r}{\ln m_-}} r^d \dpy\frac{1}{q(k)} \quad .
\eenn
Now using equation (\ref{eq=findr}), one easily checks that 
\benn
m^{k+\frac{\ln r}{\ln m_-}} \leq m^{k-k+\frac{\ln p-\ln \epsilon +\ln C_2}{\ln m_-}}&=&
C_3 p^{\frac{\ln m}{\ln m_-}} \quad .\\[0.3cm]
\lambda_p^{k+\frac{\ln r}{\ln m_-}} \geq \lambda_p^{k-k+\frac{\ln p-\ln \epsilon +\ln C_1}{\ln m_-}} &=& C_4 p^{\frac{\ln \lambda_p}{\ln m_-}} \quad .
\eenn
So
\bsenn
\dpy\frac{\nu(B_r(x))}{r^{\beta-d}} \leq K_2 p^{\frac{\ln m-\ln \lambda_p}{\ln m_-}} q^{-1} r^d \quad .
\esenn
Thus the $r^d$-term dominates. 
\QED

\pagebreak

\section{The Hausdorff Outer Measure}
\label{s-sets}
\setcounter{figure}{0}
\setcounter{equation}{0}
 
\vskip.2in

In this section we give some results about the $\beta$-dimensional Hausdorff 
outer measure of $\delta\Lambda \cap V$. These results are essentially corollaries 
of the calculations done before. We will deal only with the case where the 
eigenvalues of $M$ have equal modulus. Thus 
\bse
m_-^n=m
\label{eq=equalevalues}
\ese

\noindent
{\bf Proof of theorem \ref{thm=ssets}:} From the proof of proposition \ref{thm=upper} and lemma \ref{thm=countdelta}, we have that 
\benn
{\cal H}_{\epsilon_k}^{\beta} (\delta\Lambda \cap V) &\leq& C p^{\beta}m_-^{-k\beta} q \lambda_p^k \\[0.3cm]
&=& C p^{\beta}q \quad .
\eenn
Thus
\bsenn
{\cal H}^\beta (\delta\Lambda \cap V) \leq C \quad .
\esenn

For the second and third statements, we see that the proof of theorem \ref{thm=lowerbound}
together with equalities (\ref{eq=equalevalues}) and (\ref{eq=defbeta}) imply 
that
\bsenn
{\cal H}^\beta (\delta\Lambda \cap V) \geq \dpy\frac{\nu(\delta\Lambda \cap V)}{p^{n-\beta}q^{-1}} \quad ,
\esenn
and recall that $p$ has degree $d_M-1$ and $q$ has degree $d_{\lambda_p}-1$.
\QED

\pagebreak

\section{Some Examples}
\label{examples}
\setcounter{figure}{0}
\setcounter{equation}{0}
 
\vskip.2in

We illustrate the ideas in this work by calculating the dimension of the boundary 
of a self-affine tile in a number of examples.

We begin with an trivial example that can be understood without calculation. 

\begin{exam}
\label{example1}
The interval $[0,1]$ is the invariant set for $(M,R)= (2,\{0,1\})$. It is easy 
to check that D=\{-1,0,1\} with $0$ appearing with multiplicity 2. The set 
$S$ is given by $\{-1,0,1\}$ and $S^+= \{0,1\}$. Thus the 
transition matrix and the reduced transition matrix are given by 
\bsenn
T = \left( \begin{array}{lll} 1&1&0\\0&2&0\\0&1&1 \end{array} \right) \quad 
; \quad T^+ = \left( \begin{array}{ll} 2&0\\1&1 \end{array} \right) \quad .
\esenn
We follow the convention that the ordering of coordinates on which $T^+$ ($T$) acts is 
the same as the ordering in $S^+$ ($S$). So the upper right entry of the $T^+$ 
multiplies $v_{\{0\},\{1\}}$, and so forth. 
Since the only special eigenvalue (definition \ref{def=specialevalues})  
equals 1, the dimension of the boundary equals 
$\dpy\frac{\ln 1}{\ln 2}=0$.
\end{exam}

Now let us look at some non-trivial examples. 

\begin{exam}
\label{exaqmple2}
Let $(M,R)= (3,\{0,4,11\})$. Then $\Hdim \delta\Lambda \cap V \approx \dpy\frac{\ln 2.84\cdots}{\ln 3} \approx 0.87\cdots$ . 
\end{exam}

\noindent
{\bf Proof:} Check that 
\bsenn
S^+ = \{0,1,2,3,4,5\} \quad .
\esenn
With the same convention as before (noting that $D$ consists of the numbers 
-11, -7, -4, 0, 4, 7, and 11, and 0 appearing with multiplicity 3),
\bsenn
T^+ = \left( \begin{array}{llllll}
3&0&0&0&0&0\\0&1&1&0&1&0\\0&0&1&2&0&0\\0&3&0&0&0&0\\1&1&0&0&0&1\\0&0&1&1&1&0
\end{array} \right) \quad .
\esenn
Using Maple, we find that the only special eigenvalue is approximately 2.84. 
\QED

We now exhibit examples in one dimension whose boundaries have dimension 
approximating (but not equal to) 0. 

\begin{exam}
\label{example3}
For $m\geq 4$, let $(M,R)=(M,\{0,2,3,\cdots m-1,m+1\})$. Then $\Hdim (\delta\Lambda \cap V) = 
\dpy\frac{\ln 3}{\ln m}$. 
\end{exam}

\noindent
{\bf Proof:} We easily see that $S^+=\{0,1\}$. For the set $D$ the following 
holds
\benn
0 &\in D \with \multip& m\\[.3cm]
1 & "& m-3\\[.3cm]
m-1 & " & 1\\[0.3cm]
m+1 & " & 1
\eenn
Thus
\bsenn
T^+ = \left( \begin{array}{ll} m&0\\m-3&3 \end{array} \right) \quad .
\esenn
\QED

Here is a family of examples such that their boundaries have dimension 
converging to (without being equal to) 1. 

\begin{exam}
\label{example4}
Let $m\geq 4$ be even. Let 
\bsenn
(M,R) = (m,\{0,2,4,\cdots (m-4), (m-2), (m+1), (m+3), \cdots (2m-1)\}) \quad . 
\esenn
Then
\bsenn
\Hdim (\delta\Lambda \cap V) = \dpy\frac{\ln [(m-1)+((m-1)^2 + 8)^{1/2}] - \ln 2}{\ln m} \quad .
\esenn
\end{exam}

\noindent
{\bf Proof:} Again, it is easy to see that $S^+=\{0,1,2\}$. Concerning $D$, we 
only need the following information:
\benn
0 &\in D \with \multip& m\\[.3cm]
2 & " & m-2 \\[.3cm]
m-2 & " & 2\\[.3cm]
m-1 & " & m/2-1 \\[.3cm]
m+1 & " & m/2\\[.3cm]
m+2 & " & 0
\eenn
Thus
\bsenn
\T^+ = \left( \begin{array}{lll} m&0&0\\0&m-1&1\\m-2&2&0 \end{array} \right) \quad .
\esenn
The only special eigenvalue is $\dpy\frac{(m-1)+((m-1)^2 + 8)^{1/2}}{2}$.
\QED

We have not been able to find any one-dimensional examples in which the leading 
special eigenvalue is associated with a Jordan block.

In studying the dimension of boundaries of tiles in two (or more) dimensions, 
three new aspects arise. In the first place, not all expanding map are 
(conjugate to) similarities. So we may obtain estimates rather than equalities. 
That this is inevitable is clear from the following example slightly modified
from \cite{Fa2}.

\begin{exam}
Let $N>2$ and 
\bsenn
(M,R) = \left( \left( \begin{array}{ll} 2&0\\0&N \end{array} \right),
\{(0,0),(-2\lambda,N-1)\} \right) \quad .
\esenn
We calculate the dimension of $\Lambda$ (not its boundary). When $\lambda=0$, 
this gives $\Hdim(\Lambda) = \dpy\frac{\ln 2}{\ln N}$. When $\lambda\neq 0$, 
we have $\Hdim(\Lambda) = 1$.
\end{exam}

\noindent
{\bf Proof:} See example 9.10 in \cite{Fa2}. 
\QED

A second problem is that, although the algorithm that determines the set $S^+$ 
terminates after a finite number (but not a priori bounded) number of steps, this 
calculation isn't nearly as straightforward as in the one-dimensional case. In fact,
to check that a given set is indeed $S^+$ is an elementary but very longwinded 
calculation, which we leave to reader in the examples below. In 
the following examples, all eigenvalues of $M$ have equal modulus.

\begin{exam}
\label{example6}
The case where $m=2$ in two dimensions. As explained in \cite{HSV}, there are 
(modulo affine coordinate transformations) only six cases. The following three
 are representative.
\bsenn
i)\quad M=\left( \begin{array}{ll} 0&-2\\1&0 \end{array} \right), 
\quad ii)\quad M = \left( \begin{array}{ll} 0&-2\\1&1 \end{array} \right),
\quad iii) \quad M = \left( \begin{array}{ll} 0&-2\\1&2 \end{array} \right)\quad ,
\esenn
and in all three cases $R=\{(0,0),(1,0)\}$.\\
We have that in the three respective cases, the Hausdorff dimension assumes the 
values $1$, $\dpy\frac{ln 1.5\cdots}{\ln 2}$, and $\dpy\frac{ln 1.7\cdots}{\ln 2}$.
\end{exam}

\noindent
{\bf Proof:} In the three cases $D= \{(0,0), (-1,0), (1,0)\}$, the point $(0,0$ 
having multiplicity 2. Note that $(M,R)$ generates the same set 
$\Lambda$ as $(M^2,MR+R)$. In the first case, we obtain the system
\bsenn
\left( \left( \begin{array}{ll} -2&0\\0&-2 \end{array} \right), 
\{((0,0)\cup(1,0))+((0,0)\cup(0,1))\} \right) \quad . 
\esenn
It is easy to see that the 
associated set $\Lambda$ is, in fact, the unit square (by explicit substitution,
for example).

In the second case:
\bsenn
S^+ = \{(0,0),(1,0),(0,1),(1,-1)\} \quad ,
\esenn
and
\bsenn
T^+ = \left( \begin{array}{llll} 2&0&0&0\\1&0&0&1\\0&2&0&0\\0&1&1&0 \end{array} \right) \quad .
\esenn
The characteristic polynomial is $-\lambda^3 + \lambda +2$, whose only zero is
$\lambda \approx 1.5\cdots$.

In the third case: 
\bsenn 
S^+ = \{(0,0),(1,0),(1,-1),(1,-2)\} \quad , 
\esenn 
and 
\bsenn 
T^+ = \left( \begin{array}{llll} 2&0&0&0\\1&0&0&1\\0&1&1&0\\0&0&2&0 \end{array} \right) \quad . 
\esenn 
The characteristic polynomial is $-\lambda^3 + \lambda^2 +2$, whose only zero is 
$\lambda \approx 1.7\cdots$. 
\QED

Finally, we conclude by calculating the dimension of the boundaries of the two 
tiles depicted in figures 2 and 3 of \cite{HSV}. Here, the aspect arises that $M$ 
may also have Jordan blocks.

\begin{exam}
\label{example7}
Let 
\bsenn
(M,R) = \left( \left( \begin{array}{ll} 2&1\\0&2 \end{array} \right),
\{(0,0),(1,0),(0,1),(1,1)\} \right) \quad .
\esenn
Then $\Hdim(\delta \Lambda \cap V) = 1$.
\end{exam}

\noindent 
{\bf Proof:} One checks that 
\bsenn 
S^+ = \{(0,0),(1,0),(0,1),(1,-1)\} \quad ,
\esenn
and 
\bsenn
T^+ = \left( \begin{array}{llll} 4&0&0&0\\2&2&0&0\\2&0&1&1\\1&1&0&2 \end{array} \right) \quad .
\esenn
\QED

\noindent
{\bf Remark:} Notice that both $M$ and $T^+$ have a Jordan block of size 2. So it 
is unclear whether $\delta\Lambda$ is an s-set. 

\begin{exam}
\label{example8}
Let 
\bsenn
\left( \left( \begin{array}{ll} 3&0\\0&3 \end{array} \right),
\{(-1,-1),(0,-1), (1,-1), (-2,0), (0,0),(2,0),(-1,1),(0,1),(1,1)\} \right) \quad .
\esenn
Then $\Hdim(\delta \Lambda \cap V) = 1$ or $\Hdim(\delta \Lambda \cap V) = 
\dpy\frac{\ln 5}{\ln 3}$. 
\end{exam}

\noindent 
{\bf Proof:} We have
\bsenn 
S^+ = \{(0,0),(1,0),(2,0),(0,1),(1,1),(1,-1)\} \quad , 
\esenn 
and  
\bsenn 
T^+ = \left( \begin{array}{llllll} 9&0&0&0&0&0\\4&5&0&0&0&0\\4&4&1&0&0&0\\2&4&0&3&0&0\\4&2&0&2&1&0\\4&2&0&2&0&1 \end{array} \right) \quad . 
\esenn
There are now two special eigenvalues, namely 3 and 5.
\QED

\noindent 
{\bf Remark:} This last result is in fact easy to verify by inspection of the figure.

\vfil \eject

\end{document}